\documentclass[10pt]{amsart}

\usepackage{amsmath}
\usepackage{amsthm}
\usepackage{amsfonts}
\usepackage{amssymb,latexsym}
\usepackage[all]{xy}
\usepackage{graphicx}
\usepackage[mathscr]{eucal}
\usepackage{verbatim}
\usepackage{hyperref}

\theoremstyle{plain}
    \newtheorem{thm}{Theorem}[section]
    \newtheorem{prop}[thm]{Proposition}
    \newtheorem{lemma}[thm]{Lemma}

\theoremstyle{definition}

\theoremstyle{remark}

\numberwithin{equation}{section}

\newcommand{\rar}{\ensuremath{\rightarrow}}
\newcommand{\lrar}{\ensuremath{\longrightarrow}}

\newcommand{\Hom}{\textup{Hom}}
\newcommand{\StMod}{\textup{StMod}}
\newcommand{\stmod}{\textup{stmod}}

\newcommand{\uHom}{\underline{\Hom}}

\newcommand{\lstk}[1]{\stackrel{#1}{\longrightarrow}}

\newcommand{\Tate}{\widehat{H}}
\newcommand{\Char}{\textup{char}}
\newcommand{\thick}{\textup{thick}}

\newcommand{\Ind}{\textup{Ind}}
\newcommand{\Res}{\textup{Res}}
\newcommand{\uInd}{\underline{\Ind}}
\newcommand{\uRes}{\underline{\Res}}
\newcommand{\up}{{\uparrow}^G}
\newcommand{\down}{{\downarrow}_H}
\newcommand{\Rar}{\Rightarrow}
\newcommand{\Lrar}{\Leftrightarrow}

\newcommand{\hrsmash}[1]{\makebox[0pt][r]{#1}}
\newcommand{\hlsmash}[1]{\makebox[0pt][l]{#1}}

\begin{document}

% Article information
\title[Generating hypothesis for the stable module category]
{The generating hypothesis for the stable module category of a $p$-group}
\date{28 November 2006}

% author one information
\author{David J. Benson}
\address{Department of Mathematical Sciences\\
University of Aberdeen \\
Meston Building \\
\break King's college, Aberdeen AB24 3UE, Scotland, UK}
\email{$\backslash$/b$\backslash$e/n$\backslash$s/o$\backslash$%
n/d$\backslash$j/$\backslash$ (without the slashes) at maths dot
abdn dot ac dot uk}
%\thanks{}

% author two information
\author{Sunil K. Chebolu}
\address{Department of Mathematics \\
University of Western Ontario \\
London, ON N6A 5B7, Canada}
\email{schebolu@uwo.ca}
%\thanks{}

% author three information
\author{J. Daniel Christensen}
\address{Department of Mathematics \\
University of Western Ontario \\
London, ON N6A 5B7, Canada}
\email{jdc@uwo.ca}
%\thanks{}

% author four information
\author{J\'{a}n Min\'{a}\v{c}}
\address{Department of Mathematics\\
University of Western Ontario\\
London, ON N6A 5B7, Canada}
\email{minac@uwo.ca}

% AMS information

% Abstract
\begin{abstract}
Freyd's generating hypothesis, interpreted in the stable module category of a
finite $p$-group $G$, is the statement that a map between finite-dimensional
$kG$-modules factors through a projective if the induced map on Tate
cohomology is trivial. We show that Freyd's generating hypothesis holds for a
non-trivial finite $p$-group $G$ if and only if $G$ is either $C_2$ or $C_3$. We also
give various conditions which are equivalent to the generating hypothesis.
\end{abstract}

\maketitle
\thispagestyle{empty}

%\tableofcontents

% Body of the paper
\section{Introduction}

The generating hypothesis (GH) is a famous conjecture in homotopy theory due
to Peter Freyd~\cite{freydGH}. It states that a map between finite spectra
that induces the zero map on stable homotopy groups is null-homotopic. If
true, the GH would reduce the study of finite spectra $X$ to the study of
their homotopy groups $\pi_*(X)$ as modules over $\pi_*(S^0)$. Therefore it
stands as one of the most important conjectures in stable homotopy theory.
This problem is notoriously hard; despite serious efforts of homotopy
theorists over the last 40 years, the conjecture remains open,
see~\cite{ethanGH, ethan2GH}.
%In order to gain some insight into this deep
%problem, it is therefore natural to examine its analogues in some algebraic
%settings that are formally very similar to the stable homotopy category of
%spectra, such as the derived category of a commutative ring and the stable
%module category of a finite group. The GH in the derived category $D(R)$ of a
%ring $R$ is the statement that a map
%between perfect complexes in $D(R)$ that induces the zero map on homology is chain-homotopic
%to the zero map.
Keir Lockridge~\cite{keir} showed that the analogue of the GH holds in the
derived category of a commutative ring $R$ if and only if $R$ is a von Neumann
regular ring (a ring over which every $R$-module is flat). More recently,
Hovey, Lockridge and Puninski have generalised this result to arbitrary
rings~\cite{GH-D(R)}. Lockridge's result~\cite{keir} applies to any tensor
triangulated category where the graded ring of self maps of the unit object is
graded commutative and is concentrated in even degrees. Note that this
condition is not satisfied by the stable homotopy category of spectra. So in
order to better understand the GH for spectra,
% where the homotopy groups
%$\pi_*(S^0)$ of the sphere spectrum are nonzero in both even and odd degrees,
%it is important to investigate the GH in an algebraic example where the ring
%of self maps of the unit object are non-zero in odd degrees as well. An
%important and interesting example which has this property is the stable module
we formulate and solve the analogue of Freyd's GH in the stable module
category of a finite $p$-group. Here the ring of self maps of the unit object (the trivial representation $k$)
is non-zero in both even and odd degrees.

To set the stage, let $G$ be a non-trivial finite $p$-group and let $k$ be a
field of characteristic $p$.  Consider the stable module category $\StMod(kG)$
of $G$. It is the category obtained from the category of $kG$-modules by
killing the projectives. The objects of $\StMod(kG)$ are the left
$kG$-modules, and the space of morphisms between $kG$-modules $M$ and $N$,
denoted $\uHom_{kG}(M,N)$, is the $k$-vector space of $kG$-module
homomorphisms modulo those maps that factor through a projective module.
$\StMod(kG)$ has the structure of a tensor triangulated category with the
trivial representation $k$ as the unit object and $\Omega$ as the loop
(desuspension) functor. The category $\stmod(kG)$ is defined similarly using
the finite-dimensional $kG$-modules. A key fact \cite{bencar-1992} is that the
Tate cohomology groups can be described as groups of morphisms in
$\StMod(kG)$: $\Tate^i(G,M) \cong \uHom(\Omega^i k, M)$. In this framework,
the GH for $kG$ is the statement that a map $\phi\colon M \rar N$ between
finite-dimensional $kG$-modules is trivial in $\stmod(kG)$ if the induced map
in Tate cohomology $\uHom(\Omega^i k, M) \rar \uHom(\Omega^i k, N)$ is trivial
for each $i$. Maps between $kG$-modules that are trivial in Tate cohomology
will be called \emph{ghosts}. It is shown in \cite{CCM2} that there are no
non-trivial ghosts in $\StMod(kG)$ if and only if $G$ is cyclic of order $2$
or $3$. The methods in~\cite{CCM2} do not yield ghosts in $\stmod(kG)$. In
this paper, we use induction to build ghosts in $\stmod(kG)$. Our main theorem
says:

\begin{thm}\label{thm:maintheorem}
Let $G$ be a non-trivial finite $p$-group and let $k$ be a field of
characteristic $p$. There are no non-trivial maps in $\stmod(kG)$ that are
trivial in Tate cohomology if and only if $G$ is either $C_2$ or $C_3$. In
other words, the generating hypothesis holds for $kG$ if and only if $G$ is
either $C_2$ or $C_3$.
\end{thm}

Note that the theorem implies that the GH for $p$-groups does not depend on
the ground field $k$, as long as its characteristic divides the order of $G$.

We now explain the strategy of the proof of our main theorem. We begin by
showing that whenever the GH fails for $kH$,
for $H$ a subgroup of $G$, then it also fails for $kG$. We then construct
non-trivial ghosts over cyclic groups of order bigger than $3$
and over $C_p \oplus C_p$. It can be shown easily
that the only finite $p$-groups that do not have one
of these groups as a subgroup are the cyclic groups $C_2$ and $C_3$. And for
$C_2$ and $C_3$ we show that the GH holds.

For a general finite group $G$, the GH is the statement that there are no
non-trivial ghosts in the thick subcategory generated by $k$.  When $G$ is not
a finite $p$-group, our argument does not necessarily produce ghosts in
$\thick(k)$ and the GH is an open problem.

In the last section we give conditions on a finite $p$-group equivalent to the GH.
One of them says
that the GH holds for $kG$ if and only if the category $\stmod(kG)$ is
equivalent to the full subcategory of finite coproducts of suspensions of $k$.
We also show that if the GH holds for a finite $p$-group, then the Tate cohomology
functor $\Tate^*(G, -)$ from $\stmod(kG)$ to the category of graded modules
over the ring $\Tate^*(G,k)$ is full.

Throughout we assume that the characteristic of $k$ divides the order of the
finite group $G$. For example, when we write $kC_3$, it is implicitly
assumed that the characteristic of $k$ is 3. We denote the desuspension of
$M$ in $\StMod(kG)$ by $\Omega(M)$, or by $\Omega_G(M)$ when we need to specify the group in question.
All modules are assumed to be left modules.

\section{Proof of the main theorem} \label{sec:maintheorem}

Suppose $H$ is a subgroup of $G$. A natural question is to ask how the truth
or falsity of the GH for $H$ is related to that for $G$. We begin by addressing this
question.

\begin{prop} \label{prop:liftingghosts}
Let $H$ be a subgroup of a finite group $G$. If $\phi$ is a ghost in
$\stmod(kH)$, then $\phi\up$ is ghost in $\stmod(kG)$. Moreover, if $\phi$ is
non-trivial in $\stmod(kH)$, then so is $\phi \up$ in $\stmod(kG)$.
\end{prop}

\begin{proof}
It is well known that the restriction $\uRes_H^G$ and induction $\uInd_H^G$ functors form an adjoint pair of exact
functors; see~\cite[Cor.~5.4]{landrock} for instance. Therefore, for any $kH$-module $L$, we have a natural isomorphism
\[   \uHom_{kH} ((\Omega^i_G k) \down, L) \cong  \uHom_{kG}(\Omega^i_G k, L \up). \]
But since  $(\Omega^i_G k) \down \cong  \Omega^i_H k$ in $\stmod(kH)$, the above isomorphism can be written as
\[  \uHom_{kH} (\Omega^i_H k, L) \cong  \uHom_{kG}(\Omega^i_G k, L \up).\]
The proposition now follows from the naturality of this
isomorphism. The second statement follows from the
observation that $\phi$ is a retract of $\phi \up \down$.
\end{proof}

Proposition~\ref{prop:liftingghosts} implies that if $G$ is a finite $p$-group, then the GH fails for $kG$ whenever
it fails for a subgroup of $G$.

We now state two lemmas which will be needed in proving our main theorem.

\begin{lemma} \label{lem:x-1isaghost}
Let $G$ be a finite $p$-group and let $x$ be a central element in $G$. Then for any $kG$-module $M$, the map
$x - 1 \colon M \rar M$ is a ghost.
\end{lemma}

\begin{proof}
Since $x$ is a central element, multiplication by $x - 1$ defines a $kG$-linear map.
We have to show that for all $n$ and all maps $f \colon \Omega^n k \rar M$, the
composition  $\Omega^n k \lstk{f} M \lstk{x-1} M $   factors through a
projective. To this end, consider the commutative diagram
\[
\xymatrix{
\Omega^n k  \ar[r]^f \ar[d]_{x -1}& M \ar[d]^{x -1} \\
\Omega^n k \ar[r]_f & M
 . }
\]
Note that $x -1$ acts trivially on $k$, so  functoriality of $\Omega$ shows that the left vertical map is stably trivial.
By commutativity of the square, the desired composition factors through a projective.
\end{proof}

\begin{lemma} \label{lemma:GHranktwogroups}
Let $G$ be a finite $p$-group and  let $H$ be  a non-trivial proper normal subgroup of $G$. If $x$ is a central element in
$G - H$, then multiplication by $x - 1$ on $k_H \up$ is a
non-trivial ghost, where $k_H$ is the trivial $kH$-module. In particular, the GH fails for $k(C_p \oplus C_p)$.
\end{lemma}

\begin{proof}
Since $k_H\up \down$ is a trivial $kH$-module, non-triviality of $x - 1$ is
easily seen by restricting to $H$. The fact that $x - 1$ is a
ghost follows from Lemma~\ref{lem:x-1isaghost}.
The last statement follows because $k_H\up$ is finite-dimensional.
\end{proof}

\begin{proof}[Proof of Theorem \ref{thm:maintheorem}]
If $G \cong C_2$ and $\Char\, k=2$, then $kC_2 \cong k[x]/(x^2)$,
so by the structure theorem for modules over a PID it is clear that every
$kG$-module is stably isomorphic to a sum of copies of $k$.
Similarly, if $G \cong C_3$ and $\Char\, k = 3$,
then one sees that every $kG$-module is stably isomorphic to a sum of copies of $k$ and $\Omega(k)$.
It follows that there are no non-trivial ghosts between finite-dimensional $kG$-modules if $G$
is either $C_2$ or $C_3$.

Now suppose that $G$ is not isomorphic to $C_2$ or $C_3$.
It suffices to show that in these cases the GH fails for some subgroup of $G$.
It is an easy exercise to show that if $G$ is not isomorphic to $C_2$ or $C_3$,
then  $G$ either has a cyclic subgroup of order at least four, or a subgroup isomorphic to
$C_p \oplus C_p$ for some prime $p$. In Lemma~\ref{lemma:GHranktwogroups} we have seen that the GH fails for
$k(C_p \oplus C_p)$. We will be done if we can show that the GH fails for cyclic groups of order at least $4$.

So let $G$ be a cyclic group of order at least $4$. Let $\sigma$ be a generator for $G$ and
let $M$ be a cyclic module of length two generated by $U$, so we have $(\sigma -1)^2 U = 0$.
Consider the map $h\colon M \rar M$ which multiplies by $\sigma - 1$:
\newdir{ }{{}*!/-5pt/@{}}
\newdir{> }{{}*!/5pt/@{>}}
\[
\xymatrix{
*{\hrsmash{$U\,$} \bullet} \ar@{{-}{-}{-}}[d]_{\sigma -1} \ar@{ -> }[drr]^h & & *{\hrsmash{$U\,$}
\bullet} \ar@{{-}{-}{-}}[d]^{\sigma -1} \\
*{\bullet} & & *{\bullet\hlsmash{~.}}
}
\]
It is not hard to see that $h$
is non-trivial, i.e., that it does not factor through the projective cover of
$M$; this is where we use the hypothesis $|G| \ge 4$. The fact that $h$ is a ghost follows from
Lemma \ref{lem:x-1isaghost}.
\end{proof}

\section{Conditions equivalent to the generating hypothesis}

\begin{thm} The following are equivalent statements for a non-trivial finite $p$-group $G$.
\begin{enumerate}
\item $G$ is isomorphic to $C_2$ or $C_3$.
\item There are no non-trivial ghosts in $\stmod(kG)$.  That is, the
  GH holds for $kG$.
\item There are no non-trivial ghosts in $\StMod(kG)$.
\item $\stmod(kG)$ is equivalent to the full subcategory of the
collection of finite  coproducts of suspensions of $k$.
\item $\StMod(kG)$ is equivalent to the full subcategory
of  arbitrary coproducts of suspensions of $k$.
%\item The covariant functor $\uHom_{kG}(\Omega^* \, k, -)$ from $\stmod(kG)$ to the category of graded (right)
%modules over $\Tate^*(G, k)$ is full and faithful.
\end{enumerate}
\end{thm}

\begin{proof}
We have already seen that the statements (2) and  (4) are equivalent to (1).
The implications $(5)  \Rar (3) \Rar (2)$ are obvious. So we will be done if
we can show that $(1) \Rar (5)$. This follows immediately from the following
more general fact, due to Crawley and J{\'o}nsson~\cite{crawley-jonsson},
which was also proved independently by Warfield~\cite{warfield}. It states
that if $G$ has finite representation type (i.e., the Sylow $p$-subgroups are
cyclic), then every $kG$-module is a direct sum of finite-dimensional
$kG$-modules.
% In particular,when $G$ is either $C_2$ or $C_3$, every object in $\StMod(kG)$ is isomorphic to a direct sum
%of suspensions of $k$.
\end{proof}

We now state a dual version of the previous theorem. A map $d \colon M \rar N$ between $kG$-modules is
called a \emph{dual ghost} if the induced map
\[\uHom_{kG}(M,  \Omega^i k) \longleftarrow \uHom_{kG}(N , \Omega^i k)\]
is zero for all $i$.

\begin{thm} The following are equivalent statements for a non-trivial finite $p$-group $G$.
\begin{enumerate}
\item $G$ is isomorphic to $C_2$ or $C_3$.
\item[$(2')$]  There are no non-trivial dual ghosts in $\stmod(kG)$.
\item [$(3')$] There are no non-trivial dual ghosts in $\StMod(kG)$.
\item [$(4')$] $\stmod(kG)$ is equivalent to the full subcategory of the
collection of finite  products of suspensions of $k$.
\item [$(5')$] $\StMod(kG)$ is equivalent to the full subcategory
of retracts of  arbitrary products of suspensions of $k$.
%\item [$(6')$] The contravariant functor $\uHom_{kG}(-,\Omega^*\, k)$ from $\stmod(kG)$ to the category of graded
%(left) modules over $\Tate^*(G, k)$ is full and faithful.
\end{enumerate}
\end{thm}

\begin{proof}
Every  finite-dimensional $kG$-module $M$ is naturally isomorphic to its double dual $M^{**}$. Therefore, the exact functor
$M \mapsto M^*$ gives a tensor triangulated equivalence between $\stmod(kG)$ and  its opposite category. This shows that
 $(2')\Lrar (2)$. In any additive category finite coproducts and finite products are the same, therefore
 $(4') \Lrar (4)$. Thus, statements $(1)$, $(2')$, and $(4')$ are equivalent.
We will be done if we can show that $(5') \Rar (3') \Rar (1) \Rar (5')$.

$(5') \Rar (3')$: Fix an arbitrary $kG$-module $M$. We have to show that there
are no non-trivial dual ghosts out of $M$. Consider the full subcategory of
all modules $X$ such that there is no non-trivial dual ghost from $M$ to $X$.
This subcategory clearly contains arbitrary products of  suspensions of $k$
and is closed under taking retractions. So by assumption the subcategory has
to be $\StMod(kG)$.

$(3') \Rar (1)$: $(3')$ clearly implies $(2')$. But we have already observed that $(2') \Rar (2) \Rar (1)$.

$(1) \Rar (5')$: We know that $(1) \Rar (5)$. It remains to show that $(5)
\Rar (5')$. Let $M$ be any $kG$-module. By $(5)$, $M$ is a  coproduct
$\oplus\, \Omega^s k$ of suspensions of $k$. We will complete the proof by
showing that the canonical map
\[ \bigoplus \; \Omega^s k \lrar \prod  \; \Omega^s k \]
is a split monomorphism in $\StMod(kG)$. By $(5)$, the fibre $F$ of this map
is  a coproduct $\oplus\, \Omega^t k$ of suspensions of $k$.  Since the
objects $\Omega^t k$ are compact, one can show that the map $F \rar \oplus\,
\Omega^s k$ is zero and therefore the desired splitting exists.
\end{proof}

\begin{comment}
Call a map $d \colon M \rar N$ between $kG$-modules  a \emph{dual ghost}
if the induced map
\[\uHom_{kG}(M,  \Omega^i k) \longleftarrow \uHom_{kG}(N , \Omega^i k)\]
is zero for all $i$. We claim that statements (2) and (3) of the above
theorem are also equivalent to their dual statements $(2')$ and $(3')$ which are
formulated using dual ghosts instead of ghosts. To see this, we consider
the natural isomorphism
\[ \Hom_{kG} (N, \Hom_k(T, k)) \cong \Hom_{kG}(T, \Hom_k(N, k)), \]
where $T$ is a Tate resolution of $k$ and $N$ is any $kG$-module; see~\cite[Proposition 3.1.8]{ben-1}, for instance.
Since $\Hom_k(T, k)$ is a complete
injective resolution of $k$, taking homology of the chain complexes in the
last isomorphism gives natural isomorphisms
\[ \uHom_{kG}(N, \Omega^{-n} k) \cong \uHom_{kG} (\Omega^n k, N^*). \]
This implies that a map $d \colon M \rar N$ is a dual ghost if and
only if $d^* \colon N^* \rar M^*$ is a ghost.
Since the functor $M \mapsto M^*$ is faithful, it follows that $(3) \implies (3')$.  Moreover, the functor $M \mapsto M^*$
is a natural equivalence when restricted to finite-dimensional modules. Therefore
it is also clear that $(2) \iff (2')$. The implication
$(3') \implies (2')$ is trivial, and $(2) \implies (3)$  is shown in the above theorem.
This completes the proof of our claim.
\end{comment}

We end with a final observation. In the stable homotopy category of spectra,
the GH says that the stable homotopy functor $\pi_*(-)$ from the category of
finite spectra to the category of graded modules over the homotopy ring
$\pi_*(S^0)$ of the sphere spectrum is faithful.  Freyd showed~\cite{freydGH}
that if the GH is true, then  $\pi_*(-)$ is also full.  So it is natural to
ask if the same is true in other algebraic settings in which the GH is being
studied. Very recently, Hovey, Lockridge and Puninski~\cite{GH-D(R)} have
given an example of ring $R$ for which the homology functor $H_*(-)$ from the
category of perfect complexes of $R$-modules to the category of graded
$R$-modules is faithful, but not full. It turns out that from this point of
view, the stable module category of a group behaves more like  the stable
homotopy category of spectra than the derived category of a ring. More
precisely, we have the following result.

\begin{thm} \label{thm:tatefull}
Let $G$ be a finite $p$-group and let $k$ be a field of
characteristic $p$. If the GH holds for $G$, then the
functor $\Tate^*(G, -)$ from  $\stmod(kG)$ to the
category of graded modules over the graded ring $\Tate^*(G, k)$ is full.
\end{thm}

\begin{proof}
We know by Theorem~\ref{thm:maintheorem} that $G$ has to be either $C_2$ or $C_3$.
Therefore every finite-dimensional $kG$-module
$M$ is stably isomorphic to a finite sum of suspensions of $k$. In particular,  $\Tate^*(G, M)$  is a free
$\Tate^*(G, k)$-module of finite rank. It follows that the induced map
\[ \uHom_{kG}(M, X) \lrar \Hom_{\Tate^*(G, k)} (\Tate^*(G, M), \Tate^*(G, X)) \]
is an isomorphism for all $kG$-modules $X$.  Since $M$ was an arbitrary
finite-dimensional $kG$-module, we have shown that the functor $\Tate^*(G, -)$
is full, as desired.
\end{proof}

% \bibliographystyle{plain}
% \bibliography{../lit}

\end{document}